\documentclass[3p]{elsarticle}

\usepackage{amsmath}
\usepackage{amsthm}
\usepackage{amssymb}
\usepackage{bm}
\usepackage{textcomp}
\usepackage{accents}
\usepackage{color}
\usepackage{nicematrix}
\usepackage{booktabs}
\usepackage[pdftex,plainpages=false,colorlinks=true,citecolor=RoyalBlue,linkcolor=Red,urlcolor=black,filecolor=black,bookmarksopen=true]{hyperref}

\newtheorem{thm}{Theorem}

\newdefinition{rmk}{Remark}

\graphicspath{{figs/}}

\begin{document}

\begin{frontmatter}

\title{Numerical solution of the Cauchy problem for second-order evolutionary equations with damped sinusoidal memory\tnoteref{label1}}

\tnotetext[label1]{The work was supported by the Russian Science Foundation (grant No. 24-11-00058).}
\author{P.N. Vabishchevich}
\ead{vab@cs.msu.ru}
\address{Lomonosov Moscow State University, 1, building 52, Leninskie Gory, 119991 Moscow, Russia}
\address{North Caucasus Federal University, 1, Pushkin str., 355017 Stavropol, Russia}

\begin{abstract}
We consider the Cauchy problem for a second-order evolutionary equation with an integral term containing the convolution of a memory kernel with the second derivative of the solution. The kernel is specified by a finite sum of damped sinusoidal functions with positive damping coefficients, amplitudes, and frequencies. A transformation of the original nonlocal problem into an equivalent local system of coupled second-order differential equations with additional variables is proposed. For the local system, an energy identity with a dissipative term is constructed and a priori estimates of the solution are obtained. On the basis of the local formulation, a three-level difference scheme with a weight parameter is considered. The decrease of the discrete energy under standard constraints on the weight parameter is proved, which ensures unconditional stability of the scheme in the corresponding energy norm. The theoretical consideration is supplemented by the results of numerical solution of a problem with memory for a one-dimensional wave equation.
\end{abstract}


\end{frontmatter}

\section{Introduction}

Evolutionary equations with memory arise in models of continuum mechanics, viscoelasticity, wave propagation, and other processes in which the current state depends on the preceding history. Their mathematical feature is the presence of Volterra integral operators relating the solution or its derivatives at the current time to values on the preceding time interval. General issues of the theory of such equations are presented in the monographs \cite{gripenberg1990volterra,pruss1993evolutionary}. Energy estimates and decay properties of solutions of second-order equations with memory were studied, in particular, in works on wave and viscoelastic models \cite{fabrizio1992mathematical,alabau2008decay}.

The numerical solution of equations with memory is complicated by the fact that direct approximation of the integral term requires taking into account the values of the solution at all previous time levels. In a direct implementation, the computational cost of processing the convolution grows with the number of time steps. One approach to efficient computation of convolutions is based on convolution quadrature and discrete operational calculus \cite{lubich1988convolution,lubich1996nonsmooth}.
Another approach is based on accelerating the computation of convolutions using special algorithms that reduce
the computational complexity and the amount of stored history
\cite{lubich2008adaptive}. General issues of the numerical solution of
Volterra equations and integro-differential equations
are also discussed in \cite{linz1985analytical,brunner1982survey}.

For kernels representable as a finite sum of exponentials, the integral term can be replaced by a system of local differential equations for additional variables. Such an internal-variable formulation is used in computational models of viscoelasticity \cite{jang2022finite}. Methods for localization and stable approximation of evolutionary equations with difference kernels were developed in \cite{vabishchevich2022volterra}. Similar ideas are used for first-order equations with memory in the derivative of the solution \cite{vabishchevich2023firstorder}. For second-order evolutionary equations, a similar approach was used in \cite{vabishchevich2022numerical}. Representations of kernels as sums of exponential functions are also used in modern numerical models of wave propagation in viscoelastic media \cite{guo2024compressing}.

In the present work, we consider more general memory kernels consisting of damped sinusoidal components. This form makes it possible to describe oscillatory memory with a finite relaxation time of each component:
\begin{equation}\label{eq:1}
k(t)=\sum_{i=1}^{m}a_i e^{-\beta_i t}\sin(\omega_i t),
\qquad a_i>0,\quad \beta_i>0,\quad \omega_i>0.
\end{equation}
Unlike positive monotone kernels, this kernel generally changes sign. Therefore, well-developed approaches based on kernel positivity are not directly applicable.

The aim of the work is to construct and analyze a local computational algorithm for second-order evolutionary equations with damped sinusoidal memory. The original problem is transformed into an equivalent system of local second-order equations. For this system, an energy identity with dissipation is established, and then a three-level difference scheme is constructed and its energy stability is proved.

The paper is organized as follows. Section 2 states the problem under consideration in its nonlocal formulation and introduces the main notation. Section 3 describes the reduction of the nonlocal problem to a local one, which is the key step for the subsequent numerical treatment. Section 4 presents the time approximation used for the resulting local problem and discusses its properties.
More general problems with memory are briefly discussed in Section 5. Section 6 reports numerical examples for the one-dimensional case, illustrating the behaviour of the proposed approach. Section 7 summarizes the obtained results and outlines directions for further work.

\section{Problem statement}

We consider the Cauchy problem for a second-order evolutionary equation
with an integral term describing a nonlocal-in-time
interaction, when
\begin{equation} \label{eq:2}
  \frac{d^2 u}{d t^2}
  +\int_0^t k(t-s)\frac{d^2 u}{d s^2}(s)\,ds
  + Au(t) = 0,
  \quad t\in(0,T].
\end{equation}
Here $u(t)$ takes values in a finite-dimensional real Hilbert space $H$ with the scalar product $(\cdot,\cdot)$
and norm $\|\cdot\|$. By $A$ we denote a linear constant (independent of time) positive self-adjoint operator acting in the space
$H$.
Equation \eqref{eq:2} is supplemented by the initial conditions
\begin{equation} \label{eq:3}
  u(0)=u_0,
  \quad \frac{d u}{d t}(0)=0,
\end{equation}
where $u_0$ is a given element of the space $H$.
The transition to more general problems with inhomogeneous right-hand sides and a second boundary condition is not fundamental, but is associated with technically complicated details that are insignificant in the context of our consideration.

In the approximate solution of initial-boundary value problems for equations with memory, the representation of the kernel as a finite sum of exponential
functions is widely used:
\begin{equation} \label{eq:4}
  k(t)=\sum_{i=1}^{m}a_i e^{-\beta_i t},
  \quad a_i>0, \quad \beta_i>0.
\end{equation}
In the theory of evolutionary
equations with memory, kernels of positive type are considered,
for which the corresponding quadratic convolution form
is nonnegative
\cite{pruss1993evolutionary,lubich1996nonsmooth}.
The exponential structure of the kernel makes it possible to introduce auxiliary
variables satisfying local differential
equations and thereby avoid direct computation of
the integral convolution at each time step
\cite{lubich2008adaptive,vabishchevich2022volterra}.

For systems with pronounced oscillatory properties,
it is of interest to consider kernels \eqref{eq:1} that directly
describe an oscillatory component of memory.
The parameter $\beta_i$ determines the decay rate of the $i$-th component, $\omega_i$ is its frequency, and $a_i$ specifies the amplitude. The kernel satisfies the condition $k(0)=0$ and generally changes sign.
For an individual term of the memory kernel $k_i(t)=a_i e^{-\beta_i t}\sin(\omega_i t)$, the estimates
\[
|k_i(t)|\leq a_i e^{-\beta_i t},\qquad
\int_0^\infty |k_i(t)|\,dt\leq\frac{a_i}{\beta_i},
\quad i =1,2, \ldots, m.
\]
The equality
\[
\int_0^\infty k_i(t)\,dt=\frac{a_i\omega_i}{\beta_i^2+\omega_i^2},
\quad i =1,2, \ldots, m ,
\]
we cannot interpret as a measure of absolute memory intensity, since due to the sign change the integral compensation may be significant.

An important difference from positive kernels is that $k(t)$ can change sign even for $a_i,\beta_i,\omega_i>0$. Therefore, from the sign of the coefficients in the expansion one cannot conclude that $k(t)\geq0$, that the convolution is of positive type, or that the original integral operator possesses the standard passivity property. This difference is essential when comparing with results for equations with kernels of positive type \cite{lubich1996nonsmooth,mclean1993numerical,mustapha2009discontinuous,pani2002galerkin}.

Direct application of a quadrature formula to the integral term of equation \eqref{eq:1} requires accumulation of the solution history.
In the numerical solution of the Cauchy problem \eqref{eq:2}, \eqref{eq:3},
methods of reduction to local
problems are of greatest interest, making it possible to reduce computational costs.
For exponential kernels, such localization is carried out
by introducing auxiliary functions satisfying
differential equations
\cite{vabishchevich2022volterra,vabishchevich2023firstorder}.

Suppose that in a direct implementation $O(d)$ operations are required for one history value at each step, and the number of steps is $N$. Then the naive computation of the convolution has total cost of order $O(dN^2)$ and requires storage of $O(dN)$ values. For fixed $m$, the local scheme uses $O((m+1)d)$ memory and $N$ solutions of systems of size proportional to $(m+1)d$. If the solver uses the sparse structure, then the cost of a step depends on the spatial operator and the solution method, but not on the length of the time history. Consequently, the advantage of localization is most noticeable for large $N$ and a moderate number of kernel components.

The present study is aimed at developing a mathematically
justified and computationally efficient approach to solving
second-order evolutionary equations with damped sinusoidal
memory. The main idea is to transform the Cauchy problem
\eqref{eq:2}, \eqref{eq:3} with kernel \eqref{eq:1} into an equivalent local system of coupled evolutionary equations by
introducing auxiliary oscillatory variables.
For the resulting extended system, it is necessary to establish
well-posedness conditions and obtain a priori estimates ensuring
stability of the solution. On this basis, numerical approximations in time are constructed and studied, ensuring unconditional stability also at the discrete level.

\section{Reduction of the nonlocal problem to a local one}

The original problem for a kernel with damped sinusoidal memory is transformed into an equivalent problem for a system of local second-order equations by introducing auxiliary variables.

\subsection{Local problem}

To transform the original equation \eqref{eq:2}, we use the formula for differentiating a convolution.
For sufficiently smooth functions $ k(t), v(t)$, the equality
\begin{equation} \label{eq:5}
  \frac{d}{dt}\int_0^t k(t-s)v(s)\,ds =
  \int_0^t k(t-s)\frac{dv}{ds}(s)\,ds
  +k(t)v(0).
\end{equation}
holds.
Setting $v(t)=du/dt$ in \eqref{eq:5} and taking into account the second
initial condition \eqref{eq:3}, we write
the integral term in the form
\[
  \int_0^t k(t-s)\frac{d^2u}{ds^2}(s)\,ds =
  \frac{d}{dt}\int_0^t k(t-s)\frac{d u}{ds}(s)\,ds.
\]

We introduce auxiliary functions
\begin{equation}\label{eq:6}
u_i(t)=\int_0^t e^{-\beta_i(t-s)}
\sin\bigl(\omega_i(t-s)\bigr)\frac{d u}{ds}(s)\,ds,
\quad i=1,\ldots,m.
\end{equation}
Then the original equation \eqref{eq:2} takes the form
\begin{equation}\label{eq:7}
\frac{d^2u}{dt^2}+Au
+\sum_{i=1}^{m}a_i\frac{du_i}{dt}=0.
\end{equation}

By direct calculations we verify that the auxiliary functions satisfy the equations
\begin{equation}\label{eq:8}
\frac{d^2u_i}{dt^2}
+2\beta_i\frac{du_i}{dt}
+(\beta_i^2+\omega_i^2)u_i
-\omega_i\frac{du}{dt}=0,\quad i=1,\ldots,m.
\end{equation}
Homogeneous initial conditions follow from \eqref{eq:6}:
\begin{equation}\label{eq:9}
u_i(0)=0,\qquad \frac{du_i}{dt}(0)=0,
\quad i=1,\ldots,m.
\end{equation}

Problem \eqref{eq:1}--\eqref{eq:3} is equivalent to the local Cauchy problem \eqref{eq:3}, \eqref{eq:7}--\eqref{eq:9}. The converse statement follows from the uniqueness of the solution of linear equation \eqref{eq:8} for $u_i$ with given initial conditions \eqref{eq:9}, and this solution coincides with the integral representation \eqref{eq:6}.

\subsection{Energy identity}

We introduce the energy functional
\begin{equation}\label{eq:10}
\mathcal E(t)= \frac12\left\|\frac{du}{dt}\right\|^2
+\frac12(Au,u) +\frac12\sum_{i=1}^{m}\frac{a_i}{\omega_i}
\left(
\left\|\frac{du_i}{dt}\right\|^2
+(\beta_i^2+\omega_i^2)\|u_i\|^2
\right).
\end{equation}

\begin{thm}\label{t-1}
Let $A$ be a positive self-adjoint operator in $H$, and let $a_i,\beta_i,\omega_i>0, \ i=1,\ldots,m.$ Then for the solution of the local problem \eqref{eq:3}, \eqref{eq:7}--\eqref{eq:9}, the energy identity
\begin{equation}\label{eq:11}
\frac{d\mathcal E}{dt}
=-2\sum_{i=1}^{m}\frac{a_i\beta_i}{\omega_i}
\left\|\frac{du_i}{dt}\right\|^2\leq0.
\end{equation}
holds. In particular, the inequality
\begin{equation}\label{eq:12}
\mathcal E(t)\leq \mathcal E(0)
=\frac12(Au_0,u_0),\qquad 0\leq t\leq T.
\end{equation}
holds.
\end{thm}

\begin{proof}
Take the scalar product of equation \eqref{eq:7} with $du/dt$. Taking into account the self-adjointness of $A$, we obtain
\begin{equation} \label{eq:13}
  \frac{1}{2}\frac{d}{dt}
  \left (
    \left\|\frac{du}{dt}\right\|^2+(Au,u)
  \right )
  +\sum_{i=1}^{m}a_i
  \left(\frac{du_i}{dt},\frac{du}{dt}\right)=0.
\end{equation}
Taking the scalar product of the $i$-th equation \eqref{eq:8}
with $(a_i/\omega_i)\,du_i/dt$, we have
\[
  \frac{1}{2}\frac{d}{dt}
  \left (
    \frac{a_i}{\omega_i}
    \left(
      \left\|\frac{du_i}{dt}\right\|^2
      +(\beta_i^2+\omega_i^2)\|u_i\|^2
    \right)
  \right ) + 2 \frac{a_i\beta_i}{\omega_i}
        \left\|\frac{du_i}{dt}\right\|^2
  =
  a_i\left(\frac{du}{dt},\frac{du_i}{dt}\right).
\]
Summing these equalities over $i$ and adding them to equality \eqref{eq:13} for the main equation, we obtain \eqref{eq:11}. The initial conditions give $\mathcal E(0)=\frac12(Au_0,u_0)$, whence \eqref{eq:12} follows.
\end{proof}

In expanded form, estimate \eqref{eq:11} is written as
\begin{equation}\label{eq:14}
\left\|\frac{du_i}{dt}\right\|^2 +(Au(t),u(t))
+\sum_{i=1}^{m}\frac{a_i}{\omega_i}
\left( \left\|\frac{du_i}{dt}\right\|^2
      +(\beta_i^2+\omega_i^2)\|u_i\|^2
    \right)
\leq (Au_0,u_0).
\end{equation}
This identity establishes the dissipativity of the extended system. By itself, it does not mean that the energy tends to zero as $t\to\infty$: a separate analysis is required to derive asymptotic decay.

Estimates follow directly from \eqref{eq:11}:
\[
\left\|\frac{du}{dt}\right\|^2 \leq (Au_0,u_0)^{1/2},\qquad
(Au(t),u(t))\leq(Au_0,u_0),
\]
and also for each $i = 1,2,\ldots, m$,
\[
\left\|\frac{du_i}{dt}\right\|^2\leq\left(\frac{\omega_i}{a_i}(Au_0,u_0)\right)^{1/2},\quad
\|u_i(t)\|\leq\left(\frac{\omega_i}{a_i(\beta_i^2+\omega_i^2)}(Au_0,u_0)\right)^{1/2}.
\]
Integrating identity \eqref{eq:11} over time, we obtain the additional estimate
\[
2\sum_{i=1}^m\frac{a_i\beta_i}{\omega_i}
\int_0^T\left\|\frac{du_i}{dt}\right\|^2\|\,dt
=\mathcal E(0)-\mathcal E(T)\leq\frac12(Au_0,u_0).
\]
Thus, weakly damped components ($\beta_i\to0+$) can preserve energy over long times.

\subsection{Vector formulation}

It is convenient to write system \eqref{eq:7}, \eqref{eq:8} in the form of a single vector evolutionary equation of second order.
For this purpose, we introduce the vector of unknowns $\bm u=\{u,u_1,\ldots,u_m\}$ and $\bm u_0=\{u_0,0,\ldots,0\}$. Consider the direct sum of spaces $\bm H=H\oplus\cdots\oplus H$, on which the scalar product and norm are
\[
(\bm u,\bm v)=(u,v)+\sum_{i=1}^{m}(u_i,v_i).
\]
The local problem \eqref{eq:3}, \eqref{eq:7}--\eqref{eq:9} is written in the form
\begin{equation}\label{eq:15}
\bm C\frac{d^2\bm u}{dt^2}
+(\bm B+\bm D)\frac{d\bm u}{dt}
+\bm A\bm u=0,
\end{equation}
\begin{equation}\label{eq:16}
\bm u(0)=\bm u_0,\qquad \frac{d\bm u}{dt}(0)=0.
\end{equation}
The operators in equation \eqref{eq:15} are block operator matrices:
\begin{equation}\label{eq:17}
\begin{aligned}
\bm A&=\operatorname{diag}\left(A,
\frac{a_1}{\omega_1}(\beta_1^2+\omega_1^2)I,\ldots,
\frac{a_m}{\omega_m}(\beta_m^2+\omega_m^2)I\right),\\
\bm B&=\begin{pmatrix}
0&a_1I&\cdots&a_mI\\
-a_1I&0&\cdots&0\\
\vdots&\vdots&\ddots&\vdots\\
-a_mI&0&\cdots&0
\end{pmatrix},\\
\bm C&=\operatorname{diag}\left(I,\frac{a_1}{\omega_1}I,\ldots,
\frac{a_m}{\omega_m}I\right),\\
\bm D&=\operatorname{diag}\left(0,
\frac{2a_1\beta_1}{\omega_1}I,\ldots,
\frac{2a_m\beta_m}{\omega_m}I\right).
\end{aligned}
\end{equation}
For our problem \eqref{eq:1}--\eqref{eq:3}, we have
\begin{equation}\label{eq:18}
\bm A=\bm A^*>0,\quad
\bm B=-\bm B^*,
\quad \bm C=\bm C^*>0,\quad \bm D=\bm D^*\geq0.
\end{equation}
Taking the scalar product of \eqref{eq:15} with $d \bm u /d t$ and taking into account the skew-symmetry of $\bm B$, we obtain
\begin{equation}\label{eq:19}
\frac{d}{dt}\left (
\frac12 \left (\bm C\frac{d\bm u}{dt},\frac{d\bm u}{dt}\right )+
\frac12(\bm A\bm u,\bm u)
\right ) =- \left (\bm D\frac{d\bm u}{dt},\frac{d\bm u}{dt} \right )\leq0.
\end{equation}
In component form, equality \eqref{eq:19} with \eqref{eq:17} taken into account gives \eqref{eq:11}.

\section{Time approximation}

The key point in justifying the computational algorithm is related to establishing the stability of the approximate solution and an efficient computational implementation of the transition to a new time level.

\subsection{Difference scheme for the vector problem}

We introduce a uniform time grid with step $\tau$ and let
\[
y^n=y(t^n), \quad t^n = n\tau, \quad n =0, 1, \ldots , N,
\quad n \tau = T.
\]
To approximate the Cauchy problem \eqref{eq:2}, \eqref{eq:15}, \eqref{eq:16}, we use a three-level difference scheme with weight $\sigma = \mathrm{const}  > 0$:
\begin{equation}\label{eq:20}
 \bm C \frac{\bm y^{n+1} - 2 \bm y^{n} + \bm y^{n-1}}{\tau^2 }
 + (\bm B + \bm D) \frac{\bm y^{n+1} - \bm y^{n-1}}{2\tau } + \bm A \bm y^{n+\sigma} =
 \bm f^{n},
 \quad n = 1,2,\ldots,N-1,
\end{equation}
\begin{equation}\label{eq:21}
 \bm y^0 = \bm u_0 ,
 \quad \bm y^1 = \accentset{\sim }{\bm u}_1 ,
\end{equation}
using the notation
\[
 \bm y^{n+\sigma} = \sigma \bm y^{n+1} + (1-2\sigma ) \bm y^{n} + \sigma \bm y^{n-1},
 \quad   \bm y^{n} = \{y^{n}, y_1^{n}, \ldots, y_m^{n}\} .
\]
In scheme \eqref{eq:20}, \eqref{eq:21}, the second initial condition should be precomputed, for example, as follows:
\[
 \frac{y^1 - y^0}{\tau^2} +  A \frac{y^1 + y^0}{4} = 0,
\]
\[
 \frac{y^1_i - y^0_i}{\tau^2 } + \omega_i^2 \frac{y^1_i + y^0_i}{4} = 0 ,
 \quad i = 1,2, \ldots, m ,
\]
The difference scheme \eqref{eq:20}, \eqref{eq:21} approximates \eqref{eq:15}, \eqref{eq:16}
with second order in $\tau$ under sufficient smoothness of the solution $\bm u(t).$
The study of the difference problem is based on general results of the theory of stability (well-posedness) of operator-difference schemes by A.A. Samarskii \cite{SamarskiiTheory,SamarskiiMatusVabischevich2002}.

\begin{thm}\label{t-2}
Let the operators in the three-level difference scheme \eqref{eq:20}, \eqref{eq:21} satisfy conditions \eqref{eq:18}.
Then for $\sigma\geq\frac14$ and using the notation
\[
\bm w^n=\frac{\bm y^n-\bm y^{n-1}}{\tau},
\qquad
\bm z^n=\frac{\bm y^n+\bm y^{n-1}}{2}
\]
the difference solution satisfies the energy equality
\begin{equation}\label{eq:22}
E^{n+1}-E^n
=-\frac{\tau}{2}
\left(\bm D(\bm w^{n+1}+\bm w^n),
\bm w^{n+1}+\bm w^n\right)\leq0,
\end{equation}
where the discrete energy has the form
\begin{equation}\label{eq:23}
E^n= \frac 12 (\bm C\bm w^n,\bm w^n)
+ \frac 12 (\bm A\bm z^n,\bm z^n)
+ \frac 12 \left(\sigma-\frac14\right)\tau^2(\bm A\bm w^n,\bm w^n) .
\end{equation}
In particular, $E^{n+1}\leq E^n\leq E^1$, and the scheme is unconditionally stable with respect to initial data in the corresponding energy norm.
\end{thm}

\begin{proof}
Using the identity
\[
\bm y^{n+\sigma}
=\frac14(\bm y^{n+1}+2\bm y^n+\bm y^{n-1})
+\left(\sigma-\frac14\right)
(\bm y^{n+1}-2\bm y^n+\bm y^{n-1}),
\]
we rewrite the scheme \eqref{eq:20} in the form
\[
\left ( \bm C+\left(\sigma-\frac14\right)\tau^2\bm A \right )
\frac{\bm w^{n+1}-\bm w^n}{\tau}
+(\bm B+\bm D)\frac{\bm w^{n+1}+\bm w^n}{2}
+\bm A\frac{\bm z^{n+1}+\bm z^n}{2}=0.
\]
Take the scalar product of this equality with
\[
 2 (\bm z^{n+1} - \bm z^{n}) = \tau (\bm w^{n+1} + \bm w^{n}) .
\]
Taking into account the skew-symmetry of the operator $\bm B$, we obtain
\[
\begin{aligned}
\left ( \Big (\bm C \right . & + \left . \Big (\sigma - \frac{1}{4} \Big ) \tau^2 \bm A \Big ) \bm w^{n+1}, \bm w^{n+1} \right )   -
 \left ( \Big (\bm C + \Big (\sigma - \frac{1}{4} \Big ) \tau^2 \bm A \Big ) \bm w^{n}, \bm w^{n} \right ) \\
 & +
 (\bm A z^{n+1},  z^{n+1}) - (\bm A z^{n},  z^{n}) + \frac{\tau}{2}
 \left(\bm D(\bm w^{n+1}+\bm w^n),
 \bm w^{n+1}+\bm w^n\right) = 0 .
 \end{aligned}
\]
Thus we obtain \eqref{eq:22} with energy \eqref{eq:23}. Since $\bm D$ is nonnegative, the right-hand side of \eqref{eq:22} is nonpositive. Moreover, for $\sigma\geq1/4$ all terms in \eqref{eq:23} are nonnegative, which completes the proof.
\end{proof}

In the most interesting case $\sigma = \frac 14$, scheme \eqref{eq:20} gives
\begin{equation}\label{eq:24}
\frac{y^{n+1}-2y^n+y^{n-1}}{\tau^2}
+A\frac{y^{n+1}+2y^n+y^{n-1}}{4}
+\sum_{i=1}^{m}a_i
\frac{y_i^{n+1}-y_i^{n-1}}{2\tau}=0,
\end{equation}
and
\begin{equation}\label{eq:25}
\begin{aligned}
&\frac{y_i^{n+1}-2y_i^n+y_i^{n-1}}{\tau^2}
+2\beta_i\frac{y_i^{n+1}-y_i^{n-1}}{2\tau}\\
& \qquad +(\beta_i^2+\omega_i^2)
\frac{y_i^{n+1}+2y_i^n+y_i^{n-1}}{4}
-\omega_i\frac{y^{n+1}-y^{n-1}}{2\tau}=0,
\quad i=1,\ldots,m.
\end{aligned}
\end{equation}
In this case, the discrete energy is written as
\[
E^{n+1}  = \frac 12 \left\|\frac{y^{n+1}-y^n}{\tau}\right\|^2
+ \frac 12 \left(A\frac{y^{n+1}+y^n}{2},
\frac{y^{n+1}+y^n}{2}\right)  + \sum_{i=1}^{m} E^{n+1}_i ,
\]
where
\[
 E^{n+1}_i = \frac{a_i}{2 \omega_i}
 \left (
 \left\|\frac{y_i^{n+1}-y_i^n}{\tau}\right\|^2
 +(\beta_i^2+\omega_i^2)
 \left\|\frac{y_i^{n+1}+y_i^n}{2}\right\|^2
 \right ) ,
 \quad i = 1,2, \ldots, m .
\]
Energy equality \eqref{eq:22} shows that this quantity does not increase with the time level number.

\subsection{Computational implementation}

Let us present computational formulas for finding the solution at a new time level for $\sigma = \frac 14$.
First, from \eqref{eq:25} we express $y_i^{n+1}$ in terms of $y^{n+1}$.
We write an individual equation of the system in the form
\[
p_i y_i^{n+1}+q_i y_i^n+r_i y_i^{n-1} =
\omega_i\frac{y^{n+1}-y^{n-1}}{2\tau},
\]
in which
\[
p_i=\frac1{\tau^2}+\frac{\beta_i}{\tau}+\frac{\beta_i^2+\omega_i^2}{4} >0,
\quad q_i=-\frac2{\tau^2}+\frac{\beta_i^2+\omega_i^2}{2},\qquad
r_i=\frac1{\tau^2}-\frac{\beta_i}{\tau}+\frac{\beta_i^2+\omega_i^2}{4}.
\]
Therefore,
\begin{equation}\label{eq:26}
y_i^{n+1}=\gamma_i y^{n+1}+\delta_i,
\quad i = 1,2,\ldots,m,
\end{equation}
where
\[
\gamma_i=\frac{\omega_i}{2\tau p_i} ,
\quad \delta_i =
-\gamma_i y^{n-1} -\frac{q_i y_i^n+r_i y_i^{n-1}}{p_i}.
\]

Substitution of \eqref{eq:26} into \eqref{eq:24} gives a problem for finding
$y^{n+1}$:
\begin{equation}\label{eq:27}
\left( \frac1{\tau^2} I +\frac1{2\tau}\sum_{i=1}^m a_i\gamma_i
I +\frac14 A \right) y^{n+1} = f^{n+1} ,
\end{equation}
with right-hand side
\[
 \begin{aligned}
 f^{n+1}
 & = \frac{2}{\tau^2}y^n  -\frac14 A\bigl(2y^n+y^{n-1}\bigr)
 +\left( -\frac1{\tau^2}  +\sum_{i=1}^m\frac{a_i\gamma_i}{2\tau}
 \right) y^{n-1}
 \\
 & + \sum_{i=1}^m  \frac{a_i}{2\tau p_i}
 \left( q_i y_i^n+(r_i+p_i)y_i^{n-1}  \right).
 \end{aligned}
\]
After solving equation \eqref{eq:27}, the auxiliary variables are recalculated according to \eqref{eq:26}.

With this organization of computations, the computational costs for accounting for memory effects are insignificant.
The problem for determining the main component of the solution $y^{n+1}$ is not complicated, while the auxiliary components $y_i^{n+1}, \ i = 1,2,\ldots,m,$ are recalculated explicitly.

\section{More general problems}

We consider a generalization of the original problem to the case when
the operator multiplying the second derivative and the operator entering
the integral term are different from the identity operator.

\subsection{Problem formulation}

Suppose it is required to find a function $u(t)$ with values in a
finite-dimensional real Hilbert space $H$
satisfying the equation
\begin{equation}\label{eq:28}
  B\frac{d^2u}{dt^2}
  +\int_0^t k(t-s)C\frac{d^2u}{ds^2}(s)\,ds
  +Au(t)=0,
  \quad t\in(0,T],
\end{equation}
with initial conditions \eqref{eq:3}.
We assume that the operators $A$, $B$, and $C$ are linear constant positive self-adjoint operators in $H$, and the memory kernel has the form \eqref{eq:1}.

When introducing auxiliary functions according to \eqref{eq:6}, the local problem is formulated similarly to \eqref{eq:3}, \eqref{eq:7}--\eqref{eq:9}.
Instead of \eqref{eq:7}, we have the equation
\begin{equation}\label{eq:29}
  B\frac{d^2u}{dt^2}
  +C\sum_{i=1}^{m}a_i\frac{du_i}{dt}
  +Au=0,
\end{equation}
For $\mathcal E(t)$ we obtain
\begin{equation}\label{eq:30}
\mathcal E(t)= \frac12\left(B \frac{du}{dt},\frac{du}{dt}\right)
+\frac12(Au,u) +\frac12\sum_{i=1}^{m}\frac{a_i}{\omega_i}
\left(
\left(C\frac{du_i}{dt}, \frac{du_i}{dt}\right)
+(\beta_i^2+\omega_i^2)(C u_i, u_i)
\right).
\end{equation}
Similarly to Theorem~\ref{t-1}, for problem \eqref{eq:3}, \eqref{eq:8}, \eqref{eq:9}, \eqref{eq:29}, the energy identity
\[
  \frac{d\mathcal E}{dt}
  =-2\sum_{i=1}^{m}\frac{a_i\beta_i}{\omega_i}
  \left(C\frac{du_i}{dt}, \frac{du_i}{dt}\right)\leq0.
\]
and inequality \eqref{eq:12} are proved.

\subsection{Vector problem}

The problem under consideration is written in vector form \eqref{eq:15}, \eqref{eq:16}, where now
\[
\begin{aligned}
\bm A&=\operatorname{diag}\left(A,
\frac{a_1}{\omega_1}(\beta_1^2+\omega_1^2) C,\ldots,
\frac{a_m}{\omega_m}(\beta_m^2+\omega_m^2) C\right),\\
\bm B&=\begin{pmatrix}
0&a_1 C&\cdots&a_m C\\
-a_1 C&0&\cdots&0\\
\vdots&\vdots&\ddots&\vdots\\
-a_m C&0&\cdots&0
\end{pmatrix},\\
\bm C&=\operatorname{diag}\left(B,\frac{a_1}{\omega_1}C,\ldots,
\frac{a_m}{\omega_m}C\right),\\
\bm D&=\operatorname{diag}\left(0,
\frac{2a_1\beta_1}{\omega_1}C,\ldots,
\frac{2a_m\beta_m}{\omega_m}C\right).
\end{aligned}
\]
The main properties of these operators \eqref{eq:18} are preserved.
This enables us to prove, similarly to Theorem~\ref{t-2}, a statement on the stability of the three-level operator-difference scheme
\eqref{eq:20}, \eqref{eq:21}
for the approximate solution of problem \eqref{eq:3}, \eqref{eq:8}, \eqref{eq:9}, \eqref{eq:29} for $\sigma\geq\frac14$.

\section{Numerical examples}

In the calculations presented, we restricted ourselves to one-dimensional wave problems in space.
Given the current state of numerical solution of boundary value problems for elliptic equations, the transition to more complex multidimensional wave problems with memory is not of a fundamental nature.
Numerical data are presented on sufficiently fine computational grids in time and space, where we can neglect computational errors.

\subsection{One-dimensional problem}

Memory effects in the one-dimensional oscillation equation are taken into account by using the damped sinusoidal kernel \eqref{eq:1} with $m=2$. This leads us (see \eqref{eq:3}, \eqref{eq:7}--\eqref{eq:9}) to the system of equations
\[
\frac{\partial ^2 v}{\partial  t^2} - a^2 \frac{\partial^2 v}{\partial  x^2}
+ a_1 \frac{\partial v_1}{\partial  t} + a_2 \frac{\partial v_2}{\partial  t}  = 0  ,
\]
\[
\frac{\partial ^2v_i}{ \partial t^2}
+2\beta_i\frac{\partial v_i}{\partial t}
+(\beta_i^2+\omega_i^2)v_i
-\omega_i\frac{\partial v}{\partial t}=0,
\quad i = 1,2 .
\]
The boundary and initial conditions have the form
\[
v(0,t) = v(l,t) = 0 ,
\]
\[
v(x,0) = v_0(x) ,
\quad \frac{\partial v}{\partial  t} (x,0)  = 0 ,
\]
\[
v_i(x,0) = 0 ,
\quad \frac{\partial v_i}{\partial  t}(x,0)  = 0 ,
\quad i = 1,2 .
\]
Numerical illustrations are performed by tracking the dynamics of the initial solution profile $v(x,0)$.

The approximate solution of the problem is obtained using a difference approximation in space.
We assume that a uniform grid with step $h$ is introduced on the interval $[0,l]$. Let $\omega_h$ denote the set of interior nodes of the spatial grid. For grid functions $u(x) = 0, \ x \notin \omega_h$, we define the difference operator $A$ by the relation
\[
 A u = - \frac{1}{h^2} \big (u(x+h) - 2 u(x) + u(x-h) \big ),
 \quad x \in \omega_h .
\]
The scalar product in $H$ is introduced as follows:
\[
(u,v) = \sum_{x \in \omega_h} u(x) v(x) h .
\]
The resulting semi-discrete problem for the approximate solution $u(t) \approx v(\cdot,t)$ is written in the form \eqref{eq:3}, \eqref{eq:7}--\eqref{eq:9}.

\begin{figure}[ht]
\centering
\includegraphics[width=0.75\linewidth]{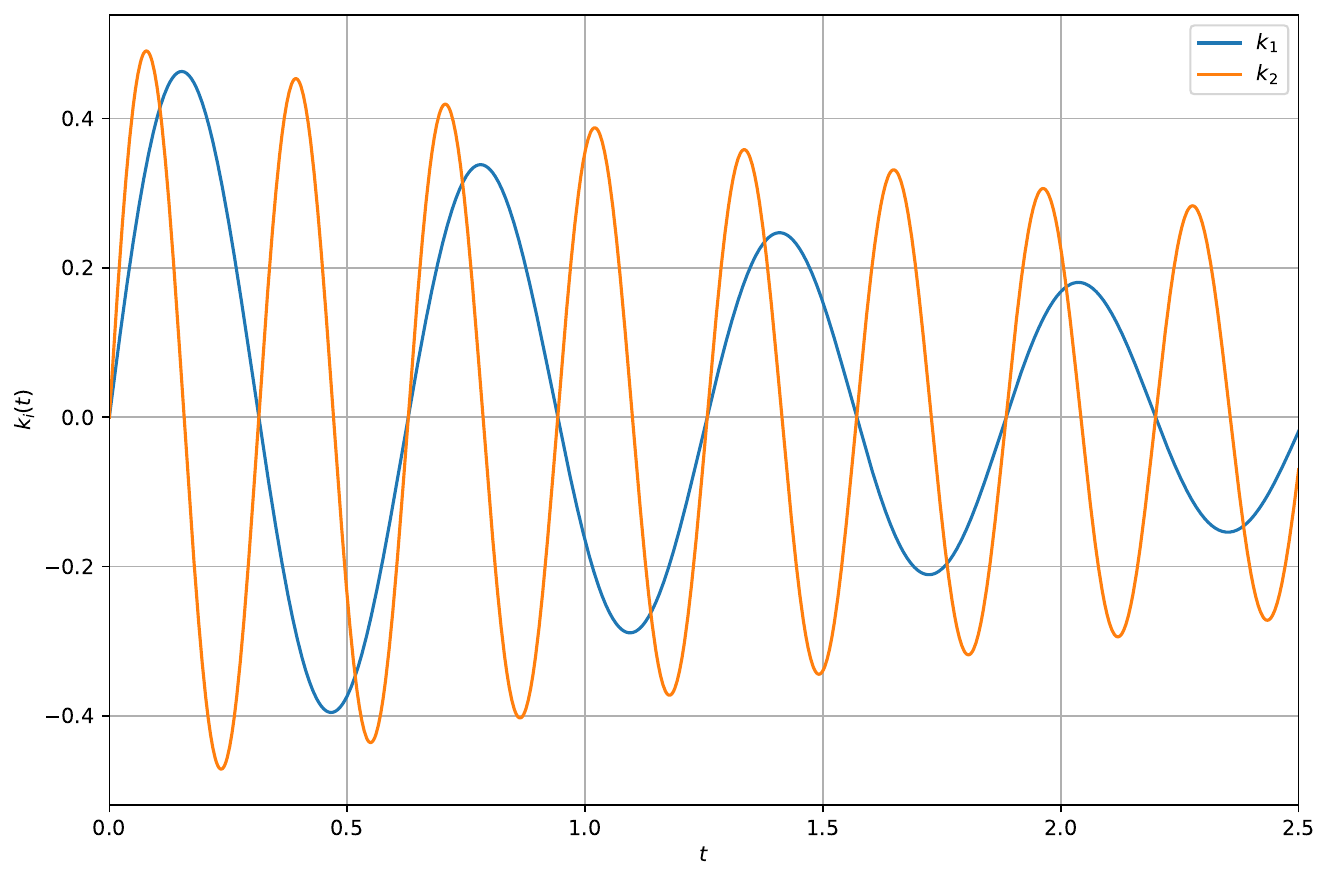}
\caption{Components of the memory kernel.}
\label{f-1}
\end{figure}

\subsection{Computational results}

The problem is considered for $l =1, \ a = 1$ and the initial condition
\[
v_0(x) = e^{- (x-x_0)^2 d^{-2}}, \quad x_0 = 0.25, \quad d=10^{-2} .
\]
Calculations were performed on a grid with $h = 10^{-3}, \ \tau = 10^{-4} .$
For the memory kernel, the representation $k(t) = k_1(t) + k_2(t)$ was used, in which
\[
\begin{split}
  k_1(t):  \  a_1 = 0.5, \ \beta_1 = 0.5, \ ~\omega_1 = 10, \\
  k_2(t):  \  a_2 = 0.5, \ \beta_2 = 0.25, \  \omega_2 = 20 .
\end{split}
\]
These dependencies are shown in Fig.~\ref{f-1}.

\begin{figure}[ht]
\centering
\includegraphics[width=0.75\linewidth]{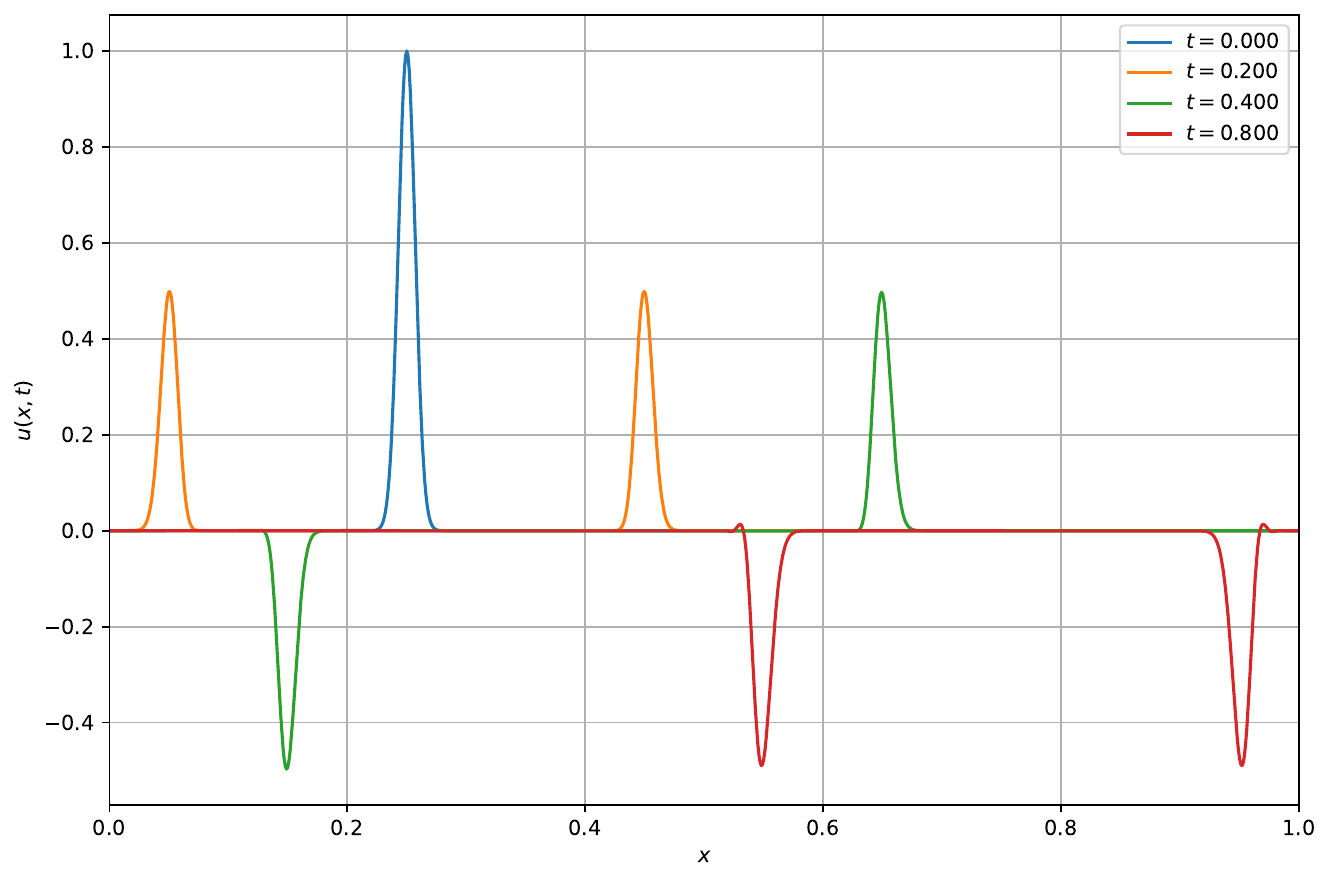}
\caption{Solution of the problem without taking memory effects into account.}
\label{f-2}
\end{figure}

\begin{figure}[ht]
\centering
\includegraphics[width=0.75\linewidth]{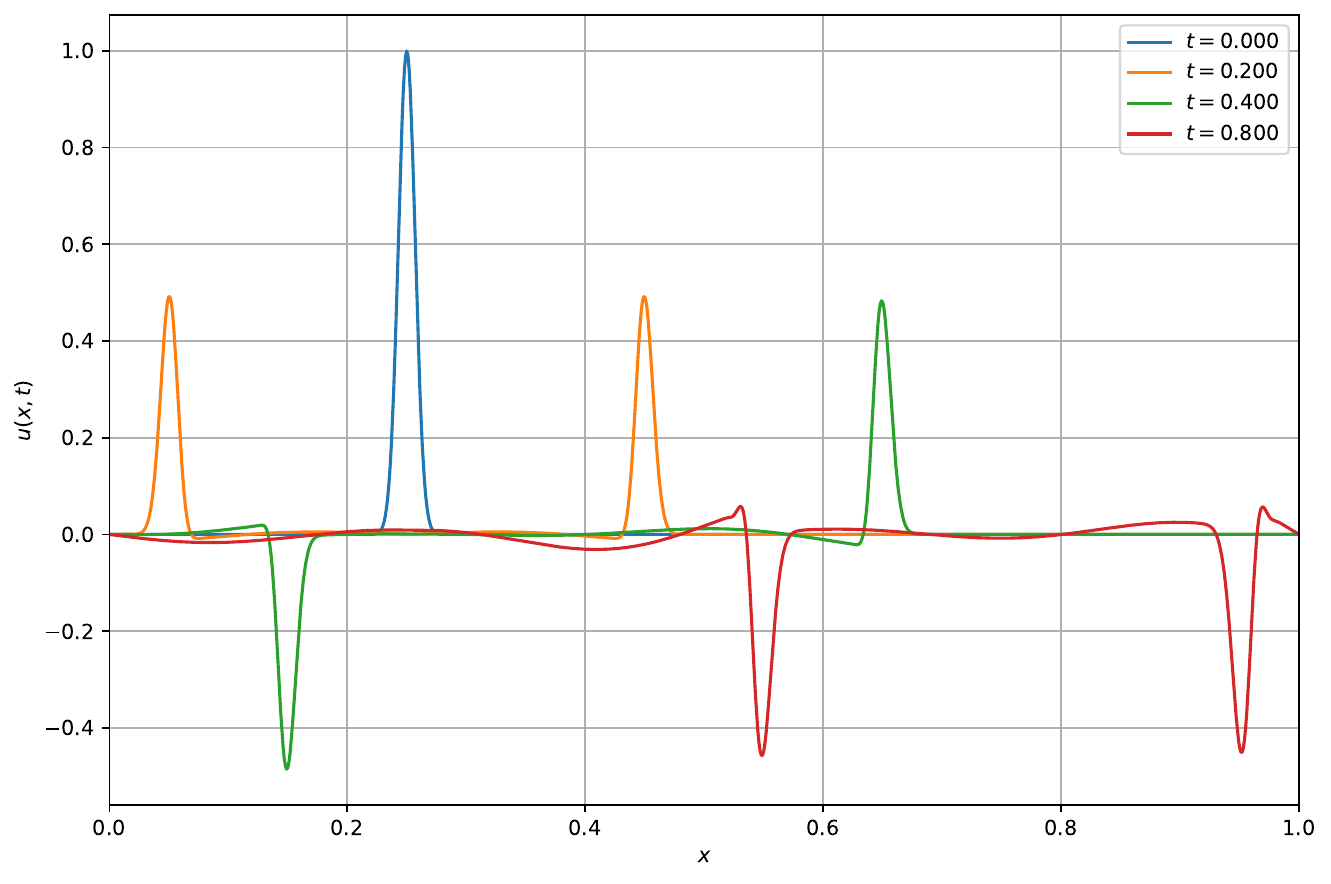}
\caption{Solution of the problem taking memory effects into account.}
\label{f-3}
\end{figure}

The numerical solution without taking memory effects into account ($a_1 = a_2 = 0$) at selected times is shown in Fig.~\ref{f-2}.
The discrete energy is conserved with high accuracy:
\[
E^n = E^1 \approx 62.51,
\quad n = 1,2, \ldots, .
\]
For the differential problem
\[
\mathcal E(0) = \frac{\sqrt{2\pi}}{4 d} \approx 62.6657 .
\]
When memory effects are taken into account, the restructuring of the solution $u(x,t)$ is illustrated in Fig.~\ref{f-3}. Comparison with Fig.~\ref{f-2} demonstrates the preservation of the dynamics of the main disturbances with a noticeable restructuring of the solution between them.

\begin{figure}[ht]
\centering
\includegraphics[width=0.75\linewidth]{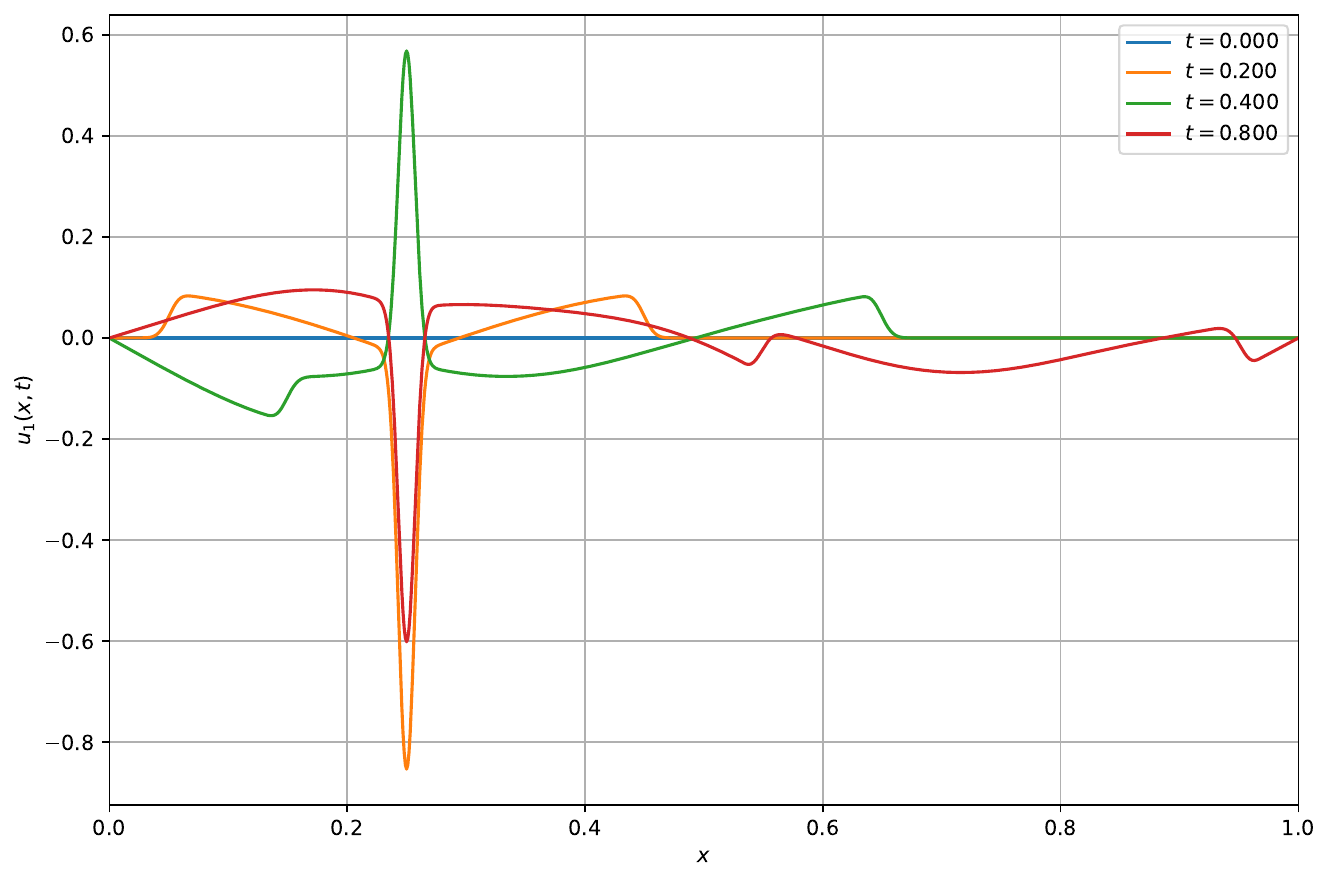}
\caption{Auxiliary function $u_1(x,t)$ at different times.}
\label{f-4}
\end{figure}

\begin{figure}[ht]
\centering
\includegraphics[width=0.75\linewidth]{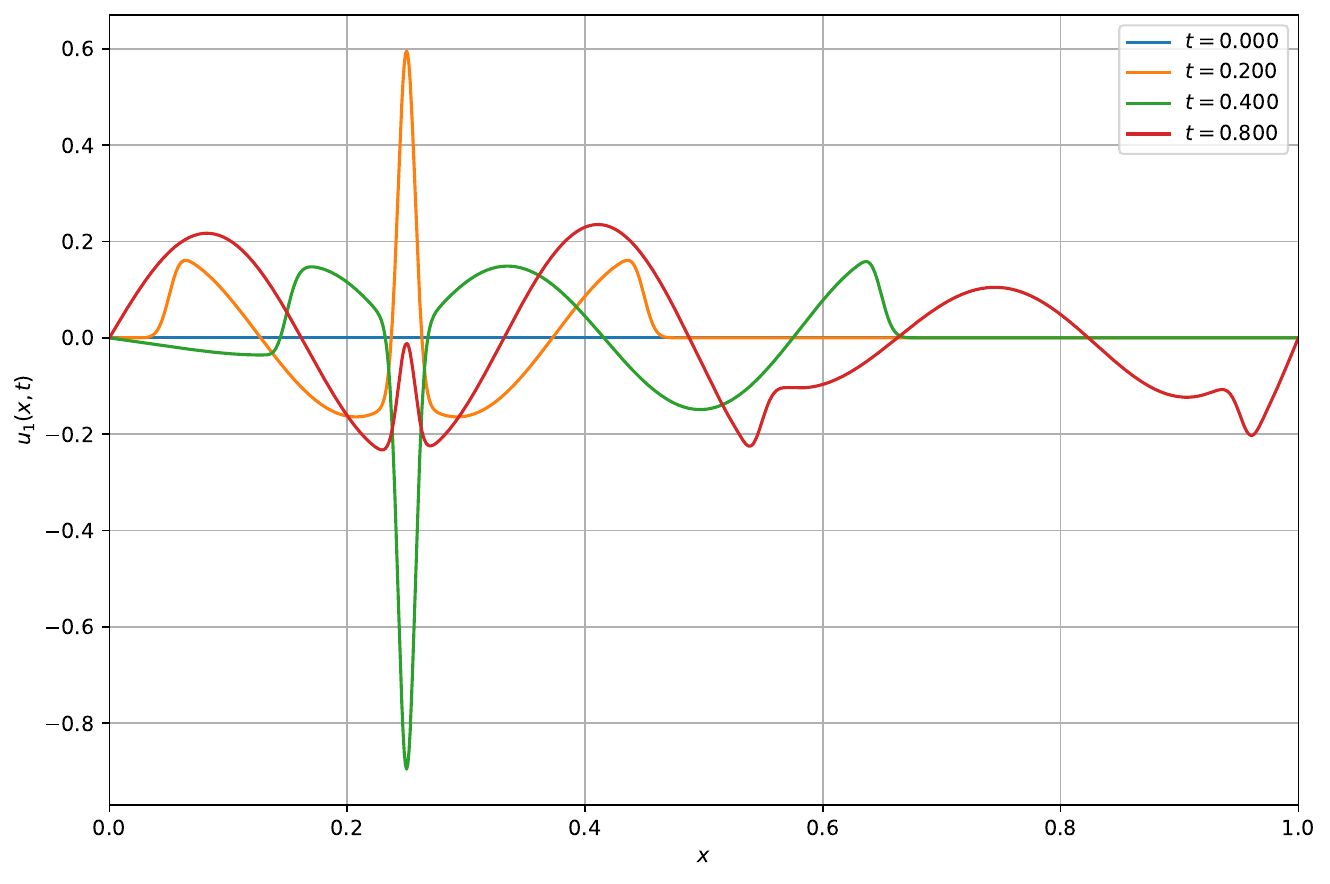}
\caption{Auxiliary function $u_2(x,t)$ at different times.}
\label{f-5}
\end{figure}

\begin{figure}[ht]
\centering
\includegraphics[width=0.75\linewidth]{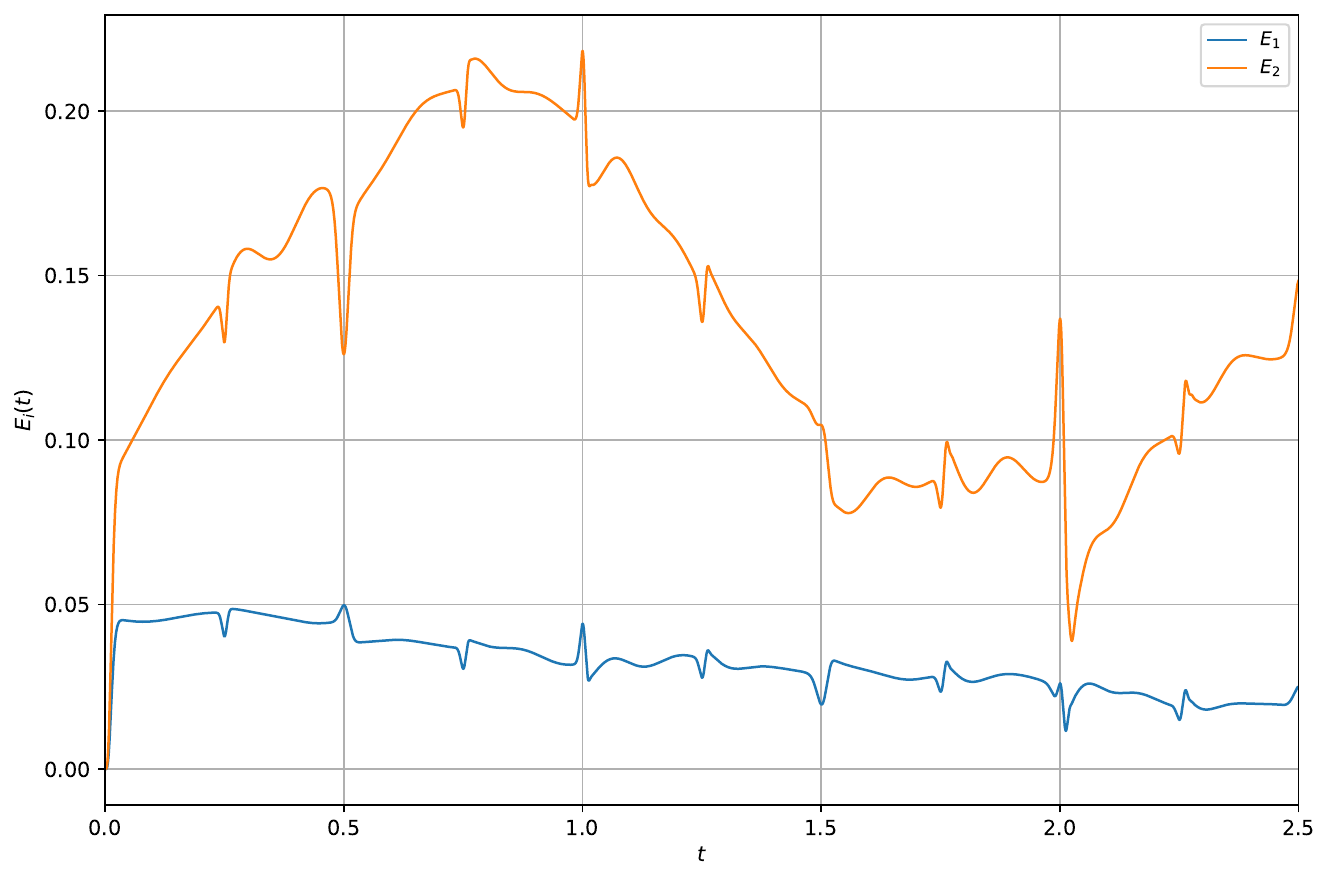}
\caption{Dynamics of individual components of the memory energy.}
\label{f-6}
\end{figure}

The dynamics of auxiliary functions $u_1(x,t)$ and $u_2(x,t)$ are shown in Fig.~\ref{f-4} and Fig.~\ref{f-5}.
They are described by equations \eqref{eq:8}.
The restructuring of these solutions is determined primarily by the first time derivative of the solution $u(x,t)$ and by the frequencies $\omega_i, \ i = 1,2$ (see the last term on the left-hand side of \eqref{eq:8}).
Therefore, oscillatory processes are less pronounced for $u_1(x,t)$ than for $u_2(x,t)$.
The contribution of individual terms $E_1(t)$ and $E_2(t)$ to the total energy of the system is illustrated in Fig.~\ref{f-6}.
An oscillatory character of the memory energy is observed, with its general decrease over time.

\clearpage

\section{Conclusion}

The work proposes and justifies a computationally efficient approach to solving second-order evolutionary equations with damped sinusoidal memory.

\begin{enumerate}[(1)]

\item The Cauchy problem for a second-order evolutionary equation with an integral term containing the convolution of a sign-changing memory kernel with the second derivative of the solution is studied. The kernel is specified by a finite sum of damped sinusoidal components, which makes it possible to describe oscillatory memory with a finite relaxation time. Unlike positive-type kernels, the kernel under consideration generally changes sign; therefore, standard analysis methods based on kernel positivity are not directly applicable.

\item The original nonlocal problem is transformed into an equivalent local system of second-order equations with auxiliary variables. For the local problem, an energy identity with a dissipative term is established and a priori estimates are obtained. The local problem is written in vector form, which made it possible, in particular, to apply the general theory of stability of operator-difference schemes.

\item For the vector problem, a three-level difference scheme with weight parameter $\sigma$ is constructed. It is proved that for $\sigma\ge 1/4$ the discrete energy does not increase, which ensures unconditional stability of the scheme with respect to initial data in the corresponding energy norm. A computational implementation is proposed in which auxiliary variables are recalculated explicitly, and the main component of the solution is determined from a problem that is not complicated compared with the problem without memory.

\item Numerical calculations for the one-dimensional wave problem with two damped sinusoidal kernel components illustrate the restructuring of the solution under the influence of memory, the dynamics of auxiliary functions, and the oscillatory character of the memory energy components with its general decrease.

\item Further research may be related to error estimates for the approximate solution, analysis of its asymptotic decay, and extension of the approach to other classes of memory kernels.
Of particular interest is work on applying the developed approach to the numerical solution of applied dynamical problems with memory.

\end{enumerate}


\begin{thebibliography}{10}
\expandafter\ifx\csname url\endcsname\relax
  \def\url#1{\texttt{#1}}\fi
\expandafter\ifx\csname urlprefix\endcsname\relax\def\urlprefix{URL }\fi
\expandafter\ifx\csname href\endcsname\relax
  \def\href#1#2{#2} \def\path#1{#1}\fi

\bibitem{gripenberg1990volterra}
G.~Gripenberg, S.-O. Londen, O.~Staffans, Volterra Integral and Functional
  Equations, Vol.~34 of Encyclopedia of Mathematics and its Applications,
  Cambridge University Press, Cambridge, 1990.

\bibitem{pruss1993evolutionary}
J.~Pr{\"u}ss, Evolutionary Integral Equations and Applications, Vol.~87 of
  Monographs in Mathematics, Birkh{\"a}user, Basel, 1993.

\bibitem{fabrizio1992mathematical}
M.~Fabrizio, A.~Morro, Mathematical Problems in Linear Viscoelasticity, SIAM,
  1992.

\bibitem{alabau2008decay}
F.~Alabau-Boussouira, P.~Cannarsa, D.~Sforza, Decay estimates for second order
  evolution equations with memory, Journal of Functional Analysis 254~(5)
  (2008) 1342--1372.

\bibitem{lubich1988convolution}
C.~Lubich, Convolution quadrature and discretized operational calculus. i,
  Numerische Mathematik 52 (1988) 129--145.

\bibitem{lubich1996nonsmooth}
C.~Lubich, I.~H. Sloan, V.~Thom{\'e}e, Nonsmooth data error estimates for
  approximations of an evolution equation with a positive-type memory term,
  Mathematics of Computation 65~(213) (1996) 1--17.

\bibitem{lubich2008adaptive}
M.~L{\'o}pez-Fern{\'a}ndez, C.~Lubich, A.~Sch{\"a}dle, Adaptive, fast, and
  oblivious convolution in evolution equations with memory, SIAM Journal on
  Scientific Computing 30~(2) (2008) 1015--1037.

\bibitem{linz1985analytical}
P.~Linz, Analytical and Numerical Methods for Volterra Equations, Vol.~7 of
  Studies in Applied Mathematics, Society for Industrial and Applied
  Mathematics, Philadelphia, 1985.

\bibitem{brunner1982survey}
H.~Brunner, A survey of recent advances in the numerical treatment of volterra
  integral and integro-differential equations, Journal of Computational and
  Applied Mathematics 8~(3) (1982) 213--229.

\bibitem{jang2022finite}
Y.~Jang, S.~Shaw, Finite element approximation and analysis of a viscoelastic
  scalar wave equation with internal variable formulations, Journal of
  Computational and Applied Mathematics 412 (2022) 114340.

\bibitem{vabishchevich2022volterra}
P.~N. Vabishchevich, Numerical solution of the cauchy problem for volterra
  integrodifferential equations with difference kernels, Applied Numerical
  Mathematics 174 (2022) 177--190.

\bibitem{vabishchevich2023firstorder}
P.~N. Vabishchevich, Approximate solution of the cauchy problem for a
  first-order integrodifferential equation with solution derivative memory,
  Journal of Computational and Applied Mathematics 422 (2023) 114887.

\bibitem{vabishchevich2022numerical}
P.~N. Vabishchevich, Numerical solution of the {C}auchy problem for a
  second-order integro-differential equation, Differential Equations 58~(7)
  (2022) 899--907.

\bibitem{guo2024compressing}
X.~Guo, S.~Jiang, Y.~Xiong, J.~Zhang, Compressing the memory variables in
  constant-{Q} viscoelastic wave propagation via an improved
  sum-of-exponentials approximation, Journal of Computational Physics 518
  (2024) 113326.

\bibitem{mclean1993numerical}
W.~McLean, V.~Thom{\'e}e, Numerical solution of an evolution equation with a
  positive type memory term, Journal of the Australian Mathematical Society.
  Series B. Applied Mathematics 35~(1) (1993) 23--70.

\bibitem{mustapha2009discontinuous}
K.~Mustapha, W.~McLean, Discontinuous galerkin method for an evolution equation
  with a memory term of positive type, Mathematics of Computation 78~(268)
  (2009) 1975--1995.

\bibitem{pani2002galerkin}
A.~K. Pani, G.~Fairweather, An $h^1$-galerkin mixed finite element method for
  an evolution equation with a positive-type memory term, SIAM Journal on
  Numerical Analysis 40~(4) (2002) 1475--1490.

\bibitem{SamarskiiTheory}
A.~A. Samarskii, The Theory of Difference Schemes, Marcel Dekker, New York,
  2001.

\bibitem{SamarskiiMatusVabischevich2002}
A.~A. Samarskii, P.~P. Matus, P.~N. Vabishchevich, Difference Schemes with
  Operator Factors, Kluwer Academic, Dordrecht, 2002.

\end{thebibliography}
\end{document}